\documentclass[11pt]{article}

 \usepackage{amssymb, amsfonts, amsmath}
 \usepackage[xdvi]{epsfig}

  \newtheorem {Theorem} {Theorem} [section]
  
  \newtheorem {Proposition} [Theorem] {Proposition}

  \newtheorem {rem} [Theorem] {}

  \newenvironment{Proof}[1][Proof] 
  {       
    \emph{#1:} 
    \setlength{\parskip}{0ex plus0ex minus0ex}}
  { \hspace*{\fill} \ensuremath{\square} \vspace{1ex}  }

    \newcommand {\norm} [2] [] {\ensuremath{ \left\Vert  #2  \right\Vert_{#1} } } 
    \newcommand {\R} {\ensuremath{\mathbb{R}}}
    
    \newcommand {\N} {\ensuremath{\mathcal{N}}}
    
    \newcommand {\F}  {\ensuremath {\mathcal{F}}}
    
    \newcommand {\om} {\omega}

    \newcommand {\T} {\mathcal{T}}
    
    \newcommand {\skalprod} [3] [] {\ensuremath{ \left\langle #2,#3 \right\rangle_{#1}}}
    \newcommand {\expval}[2] [] {\mathbb E_{#1}\left( #2 \right)}

    \newcommand {\W} {\mathcal{W}}
    
    \newcommand {\supp} {\mathrm{supp}}

    \newcommand {\mr} {\mathrm}
    \newcommand {\G}  {\ensuremath {\mathcal{G}}}
    \newcommand {\g}  {\ensuremath {\mathsf{G}}}

    \newcommand {\E} {\mathbb E}

    \renewcommand {\P} {\mathcal{P}}
    
    \newcommand {\C} {\mathbb C}
    \newcommand {\bS} {S}
   
   \newcommand{\cK} {\mathcal K}
   \newcommand{\cT} {\mathcal T}

\begin{document} 

\title
{\Large \bf A central limit theorem for Gibbs measures relative to Brownian motion }
\author{
\small Volker Betz{\thanks{Institut f\"ur Biomathematik and Biometrie, 
GSF Forschungszentrum,
Postfach 1129, D-85758 Oberschlei{\ss}heim} } and Herbert Spohn{\thanks{Zentrum Mathematik, Technische Universit\"at M\"unchen, Boltzmannstr. 3, D-85747 Garching} } \\[0.1cm]
{\small volker.betz@gsf.de, spohn@mathematik.tu-muenchen-de}
}
\date{}
\maketitle
\begin{abstract}
We study a Gibbs measure over Brownian motion with a pair potential which depends only on the increments. Assuming a particular form of this pair potential, we establish that in the infinite volume limit the Gibbs measure can be viewed as Brownian motion moving in a dynamic random environment. Thereby we are in a position use the technique of Kipnis and Varadhan and to prove a functional central limit theorem.
\end{abstract}

\section{Introduction}\label{sec.a}

We consider standard Brownian motion in $\mathbb{R}^d$, starting
at zero, weighted \`a la Gibbs as
\begin{equation}\label{a.a}
{\mathcal{N}}_{T , \mathrm{r}} = \frac{1}{Z_{T ,\mathrm{r}}}
\exp \left( -\int^T_0 \int^T_0 W (q_t - q_s\,,\, t-s)\, dt\, ds \right) \,
\W^0_\mathrm{r} \,.
\end{equation}
Here $t\to q_t$ is a Brownian path, $q_0 = 0$,  $\W^0_{\mr r}$ is the path measure of
 Brownian motion, and, in the time window $[0,T]$, the
partition function $Z_{T,\mr r}$ normalizes the weighted path measure to
one. The precise assumptions on the pair potential $W$ will be
given below, but in essence $|W(x,t)| \leq \gamma(t)$ with
$\gamma$ bounded and decaying faster than $|t|^{-3}$ at infinity.

Since in (\ref{a.a}) $W$ depends only on the increments $q_t-q_s$,
one would expect that under rescaling the weighted path measure
$\mathcal{N}_{T,\mr r}$ looks like Brownian motion with some effective
diffusion matrix $D$. To be more precise, one first has to
establish the existence of the limit measure
$$
\mathcal{N}_\mathrm{r} = \lim_{T \to \infty} \mathcal{N}_{T,
\mathrm{r}} \,.
$$
Let then $q_t$ be distributed according to $\mathcal{N}_\mathrm{r}$.
We expect the validity of the invariance principle
\begin{equation}\label{a.c}
\lim_{\varepsilon \to 0} \sqrt{\varepsilon} q_{t/\varepsilon} =
\sqrt{D} b(t)
\end{equation}
with $b(t)$ standard Brownian motion.

The conventional approach for a proof of (\ref{a.c}) is to try to
establish good mixing properties for the process of increments.
For this purpose one maps (\ref{a.a}) to a one-dimensional spin
system over $\mathbb{N}$. The single ``spin" is a continuous path
$\sigma_j (t)$, $0\leq t\leq 1$, with $\sigma_j (0)=0$. Under the
a priori measure the spins are independent and distributed
according to standard Brownian motion $\mathcal{W}^0$ over the
time span $[0,1]$. Setting $T=N$, the interaction from (\ref{a.a})
is then rewritten as
\begin{equation} \label{spin system}
\sum_{i,j=0}^{N-1} \int_{0}^1 \int_{0}^1 W (\sigma_i (t)- \sigma_j
(s)+ \Delta (i,j)\,,\, i+t-j-s)\, dt\, ds\,,
\end{equation}
where $\Delta (i,i)=0$, $\Delta (i,j)= \sum_{\ell=j}^{i-1}
\sigma_\ell (1)$ for $i>j$, and $\Delta (i,j)= \sum_{\ell=i}^{j-1}
\sigma_\ell (1)$ for $i<j$. Clearly, the path is reconstructed
from the increments as
$$
q_t=\sigma_0 (t) \quad \textrm{for} \quad 0\leq t\leq 1\,,\qquad q_t =
\sum_{j=0}^{\ell-1} \sigma_j (1) + \sigma_{\ell} (t-\ell)\quad
\textrm{for}\quad \ell \leq t < \ell +1\,,
$$
$\ell=1,\cdots, N-1$. Mixing for one-dimensional spin systems, at the level of generality needed here, is investigated by Dobrushin \cite{Do73, Do73a}. Note that the interaction in (\ref{spin system}) is many-body. Also in the applications we have in mind, $W(x,t)$ only decays like a power in the $t$-variable. Therefore it is not so obvious whether a central limit theorem for $q_t$ can be deduced with the techniques of \cite{Do73, Do73a}.

In our contribution we will prove the invariance principle by
using the Kipnis-Varadhan theorem \cite{KV}, originally developed to
deal with random motion in a random environment. This technique,
at least in its present form, requires an underlying Markov
structure. In our context it can be achieved provided $W$ has the
particular form
\begin{equation}\label{a.f}
W(x,t)= -\frac{1}{2} \int|\widehat\rho (k)|^2 e^{ik\cdot x} e^{-\omega(k)|t|}
\frac{1}{2\omega(k)} dk
\end{equation}
with 
\begin{eqnarray}
&&\omega(k)\geq 0, \quad \omega(k)=\omega(-k), \quad \widehat\rho(k)=\widehat\rho(-k)^\ast
\quad \mbox{and } \label{cond1}\\
&&\int|\widehat\rho(k)|^2 (\omega^{-1}+\omega^{-2}+\omega^{-3})dk < \infty. \label{cond2}
\end{eqnarray}
The trick is to ``linearize" the interaction in
(\ref{a.a}) by introducing the auxiliary stationary
Ornstein-Uhlenbeck process $\phi_t(x)$ with covariance
$$
{\mathbb{E}}_{\mathcal{G}} (\phi_t(x)\phi_{t'}(x')) = \int
\frac{1}{2 \omega (k)} e^{ik\cdot (x-x')}
e^{-\omega(k)|t-t'|}dk\,.
$$
Denoting its path measure by $\mathcal{G}$, (\ref{a.a}) can be
rewritten as
$$
\mathcal{N}_{T,\mathrm{r}} = \frac{1}{Z_{T,\mathrm{r}}} \E_\G \left(
\mathrm{exp} \left( - \int^T_0 \int \rho (x-q_t) \phi_t (x) \, dx \, dt \right) \right) \, \W^0_\mathrm{r} \,.
$$

According to Kipnis and Varadhan the central object is the
environment as seen from the particle, i.e. the random field
$$
\eta_t (x)= \phi_t (x + q_t)\,.
$$
In the standard applications $\eta_t$ is Markov with an
explicitly given stationary measure.  In our case, however, we
have to first take the limit $T \to \infty$ for the joint process
$ q_t$, $\phi_t (x)$. Thereby the stationary measure appears only
indirectly and is defined as the solution of an eigenvalue
problem. In Section 2, we will prove that, in the limit $T \to
\infty$, $\eta_t (x)$ is indeed a stationary, reversible Markov
process and will identify its generator. $q_t$ can be written as
an additive functional over $\eta_t (x)$ plus a martingale with
stationary increments. Following Kipnis and Varadhan, in essence,
one can thus rely on the martingale central limit theorem, with
one proviso. It must be ensured that the effective diffusion
matrix $D$ is strictly positive, see Section 4. Usually this step
requires extra considerations and is based on lower bound
estimates to the variational formula. We did not make much
progress along these lines. Instead we will rely on an idea of
Brascamp, Lebowitz and Lieb \cite{Brascamp}, 
which in the present context has been
employed before \cite{Spohn}.

We state our main theorem.
\begin{Theorem}\label{a.the1}
Define $\N_{T,\mr r}$ as in (\ref{a.a}) with $W$ given by (\ref{a.f}). 
\begin{itemize} 
\item[(i):] $\N_{T,\mr r}$ converges locally to a measure $\N_{\mr r}$ as $T \to \infty$. 
\item[(ii):] The stochastic process $q_t$ induced by $\N_{\mr r}$ satisfies a central limit theorem 
$$
\lim_{\varepsilon \to 0} \sqrt{\varepsilon} q_{t/\varepsilon} = \sqrt{D}b(t)
$$
in distribution, where $0 \leq D \leq 1$ as a $d \times d$ matrix and $b(t)$ is standard Brownian motion.
\item[(iii):] In addition to (\ref{cond1}),(\ref{cond2}) assume
$$
\int |\widehat{\varrho}(k)|^2 |k|^2 \left( \om^{-2} + \om^{-4}\right) \, dk < \infty.
$$
Then $D > 0$. 
\end{itemize}
\end{Theorem}

To make the environment process $\eta_t$ stationary, one has to
choose the symmetric time window $[-T,T]$ and to replace in
(\ref{a.a}) $\W^0_{\mathrm{r}}$ by the two-sided Brownian motion
$\W^0$ pinned at 0. This leads to the Gibbs measure
\begin{equation} \label{a.j}
\mathcal{N}_T = \frac{1}{Z_T} \exp \left( - \int^T_{-T} \int^T_{-T}
W(q_t-q_s,t-s) \, dt \, ds \right) \, \W^0\,.
\end{equation}
Let $\mathcal{N}_{T, \ell}$ be $\mathcal{N}_{T, \mathrm{r}}$
reflected at the time origin and let $\mathcal{N}_{T,\ell
\mathrm{r}}=\mathcal{N}_{T, \ell} \otimes \mathcal{N}_{T,
\mathrm{r}}$. By our assumptions on $W$, clearly,
$$
\frac{1}{c} \leq |\frac {d \mathcal{N}_{T,\ell\mathrm{r}}}{d \mathcal{N}_T}| \leq c
$$
uniformly in $T$ for some $c > 0$. Therefore a central limit theorem for
$\mathcal{N}$ is equivalent to central limit theorem for
$\mathcal{N}_\mathrm{r}$. In the sequel we will prove the
invariance principle for (\ref{a.j}). For the convenience of the
reader we restate Theorem (\ref{a.the1}) as

\begin{Theorem}\label{a.the2}
Define $\N_T$ as in (\ref{a.j}) with $W$ given by (\ref{a.f}). 
\begin{itemize} 
\item[(i):] $\N_T$ converges locally to a measure $\N$ as $T \to \infty$. 
\item[(ii):] The stochastic process $q_t$, $t \geq 0$, induced by $\N$ satisfies a central limit theorem
$$
\lim_{\varepsilon \to 0} \sqrt{\varepsilon} q_{t/\varepsilon} = \sqrt{D}b(t)
$$
in distribution, where $0 \leq D \leq 1$ as a $d \times d$ matrix, and $b(t)$ is standard Brownian motion.
\item[(iii):] In addition to (\ref{cond1}),(\ref{cond2}) assume
\begin{equation} \label{cond3}
\int |\widehat{\varrho}(k)|^2 |k|^2 \left( \om^{-2} + \om^{-4}\right) \, dk < \infty.
\end{equation}
Then $D > 0$. 
\end{itemize}

\end{Theorem}

{\bf{Remark:}} Our work is motivated by the massless Nelson model
$\cite{Nelson model, BHLMS}$. In this case $d=3$, $\omega(k)=|k|$, and $\widehat\rho$ has
a fast decay at infinity. The most stringent condition then is
$\int |\widehat\rho|^2 \omega^{-3} d^3 k < \infty$, which requires a
decay of $W$ as 
$$|W(x,t)|\leq c (1+|t|^{3+\delta})^{-1}.$$
for some $\delta>0$. Thus we cannot allow for $\widehat\rho(0)>0$
and need a mild infrared cutoff. Physically $D^{-1}$ is the
effective mass of the quantum particle when coupled to the scalar Bose
field.

\section{Diffusion driven by a stationary field}
For our particular choice of $W$, $\N$ can be written as the path measure of a diffusion driven by a stationary random field. The aim of the present section is to establish this representation.

Let $\cK_0$ be the real Hilbert space obtained by completing the subspace of $L^2(\R^d)$ on which 
\begin{equation} \label{cK}
\skalprod[\cK_0]{a}{b} = \int \widehat{a}(k) \frac{1}{2\om(k)}\widehat{b}(k)^{\ast} \, dk.
\end{equation}
is finite with respect to the inner product given by (\ref{cK}). Here $\widehat a$ denotes the Fourier transform of $a$, and $\widehat{b}(k)^{\ast}$ denotes complex conjugation of $\widehat{b}(k)$.

Let $\g$ be the Gaussian measure indexed by $\cK_0$. It will be convenient consider $\g$ on a probability space consisting of distributions. Let $A$ be a strictly positive  operator with Hilbert-Schmidt inverse in $\cK_0$, and let $\cK$ be the completion of $\cK_0$ 
with respect to the Hilbert norm 
$$ \norm[\cK]{\phi} = \norm[\cK_0]{A^{-1}\phi} \quad \forall \phi \in \cK_0.$$
Then $\g$ is supported on $\cK$, and 
$$ \E_\g (\phi(a) \phi(b)) = \skalprod[\cK_0]{a}{b}. \qquad (a,b \in D(A)),$$
where $\phi$ is considered as a linear functional defined through the limit
$$\phi(b) = \lim_{n \to \infty} \skalprod[\cK_0]{Ab}{A^{-1}\phi_n}$$ 
for any sequence $(\phi_n)$ converging to $\phi$ in $\cK$. $\cK$ is a space of distributions, depending on the choice of $A$. For our purposes, however, the special form of $A$ does not matter. 

For $a \in \cK_0$, $\phi \mapsto \phi(a)$ is an element of
$L^2(\g)$. 
Let $\G$ be the path measure of the infinite dimensional ($\cK$-valued) Ornstein-Uhlenbeck process with mean $0$ and covariance 
$$ \E_\G( \phi_{s}(a) \phi_{t}(b)) = \int \widehat{a}(k)
    \frac{1}{2\om(k)}e^{-|t-s|\om(k)} \widehat{b}(k)^\ast \, dk
    \qquad (a,b \in \cK_0).$$
$\G$ is a reversible Gaussian Markov process with reversible measure $\g$, and $t \mapsto \phi_t(a)$ is continuous for $a \in \cK_0$.

For $q \in \R^d$, let $\tau_q$ be the shift by $q$ on $\cK$, i.e. $(\tau_q\phi)(a) = \phi(a(.-q))$. More generally, for $f \in L^2(\g)$, we define $(\tau_q f)(\phi) = f(\tau_q\phi)$. $\tau_q$ is unitary on $L^2(\g)$ for each $q$, and $q \mapsto \tau_q$ is a strongly continuous group on $L^2(\g)$. 
For $T > 0$ define 
\begin{equation} \label{P-measure}
\P_T = \frac{1}{Z_T} \exp \left(- \int_{-T}^T \tau_{q_s}\phi_s( \varrho) \, ds\right) \, \W^0 \otimes \G.
\end{equation}
Let $\G^{\phi}$ denote $\G$ conditioned on $\phi_0 = \phi$, and similarly let $\W^{q}$ be two-sided Brownian motion conditioned on $q_0 = q$. $\E_{\W \otimes \G}^{q,\phi}$ denotes expectation with respect to the measure $\W^q \otimes \G^\phi$. In a similar fashion, we will use subscripts to denote the path measures and superscripts to denote conditioning throughout the paper.  
$(\ref{P-measure})$ is related to the semigroup $P_t$ given by
\begin{equation} \label{semigroup equation}
(P_t f)(q, \phi) = \E_{\W \otimes \G}^{q,\phi} \left( \exp \left(- \int_{0}^t \tau_{q_s} \phi_s( \varrho) \, ds\right) f(q_t, \phi_t) \right) 
\end{equation}
for suitable functions $f:\R^d \times \cK \to \C$. By analogy to the Feynman-Kac formula, the generator of $P_t$ can be guessed (and will turn out) to be $-H$, where
\begin{equation} \label{generator}
 Hf(q,\phi) = - \frac{1}{2} \Delta f(q,\phi) + H_{\mr f}f(q,\phi) + V_{\varrho}(q,\phi)f(q,\phi).
 \end{equation}
Above, $\Delta$ is the Laplacian on $\R^d$, $-H_{\mr f}$ is the generator of the Ornstein-Uhlenbeck process, and $V_{\varrho}(q,\phi) = \tau_q \phi(\varrho)$.

Let
$C(\R^d, L^2(\g))$ be the space of functions $f: \R^d \times \cK \to \C$ such that $q \mapsto f(q,.)$ is continuous from $\R^d$ into $L^2(\g)$, and let $C_b(\R^d,L^2(\g))$ be the subspace of functions $f \in C(\R^d,L^2(\g))$ such that
\begin{equation} \label{sup norm}
\norm[L^{\infty}(\R^d,L^2(\g))]{f} = \sup_{q \in \R^d}\norm[L^2(\g)]{f(q,.)}
\end{equation}
is finite.
We will need to study $P_t$ on two closed subspaces of $C_b(\R^d, L^2(\g))$. The first one is  
$$C_0(\R^d,L^2(\g)) = \{f \in C_b(\R^d,L^2(\g)): \lim_{|q| \to \infty} \norm[L^2(\g)]{f(q,.)} = 0 \}.$$ 
The second subspace $\cT$ is the image of $L^2(\g)$ under the operator 
$$
U:L^2(\g) \to C(\R,L^2(\g)), \quad Uf(q,\phi) = \tau_qf(\phi).
$$
Since $\tau_q$ is an isometry on $L^2(\g)$ for each $q$,
$\cT$ equipped with the scalar product
\begin{equation} \label{scalar product in T}
\skalprod[\cT]{f}{g} = \E_{\g}((U^{-1}f) (U^{-1}g)^\ast) = \skalprod[L^2(\g)]{U^{-1}f}{U^{-1}g}. 
\end{equation}
is a Hilbert space, $\norm[\cT]{f} = \norm[L^{\infty}(\R^{d},L^2(\g)]{f}$, and $U$ is an isometry from $L^2(\g)$ onto $\cT$.

\begin{Theorem} \label{strongly continuous}
$P_t$ is a strongly continuous semigroup of bounded operators on $\cT$  and on 
$C_0(\R^d, L^2(\g))$. The generator of $P_t$ on both spaces is given by $-H$, with $H$ as in (\ref{generator}).  On $\cT$, 
$H$ is a self-adjoint operator. 
\end{Theorem}
The proof is deferred to the appendix. 

\vspace{2ex}
We use $P_t$ to establish the infinite volume limits of the measures $\P_T$ and $\N_T$. To begin with, note that  
for a function $f$ that depends on $\{ q_t: -T \leq t \leq T \}$ only, 
$ \E_{\N_T}(f) = \E_{\P_T}( f).$
This can bee seen by explicitly integrating the exponential of a linear functional appearing in $\E_{\P_T}(f)$ with respect to the Gaussian measure $\G$ for fixed path $\mathsf q$. 
Let $I \subset \R$ be an interval, and let us write 
$$ \bS_I(\mathsf q)= - \int_I\int_I W(q_s-q_t, |s-t|) \, ds \, dt \quad (\mathsf q \in C(\R,\R^d))$$
in the following. We then have 
$$
 \skalprod[\cT]{1}{P_T1} = \E_\W^0 \left( e^{\bS_{[0,T]}} \right),
$$
and 
$$
\norm[\cT]{P_T1}^2 = \E_\g \left( \left( \E_{\W \otimes \G}^{0,\phi} \left( e^{-\int_0^T \tau_{q_s}\phi_s(\varrho)} \right)\right)^2 \right). 
$$
Reversing time in one of the factors inside the $\E_\g$ expectation, 
using the Markov property of $\G$ together with the fact that $\g$ is the stationary measure of $\G$ and integrating out the Gaussian field we obtain
\begin{equation} \label{gaussint3}
\norm[\cT]{P_T1}^2 = \E_\W^0( e^{\bS_{[-T,T]}}).
\end{equation}
These formulas are the key to

\begin{Theorem} \label{existence of ground state}
The infimum of the spectrum of $H$ acting in $\cT$ is an eigenvalue of multiplicity one. The corresponding eigenfunction $\Psi \in \cT$ can be chosen strictly positive.
\end{Theorem}

\begin{Proof} 
By assumption (\ref{cond2}) there exists $C_\varrho > 0$ such that
$$ \bS_{[-T,T]}(\mathsf q) \leq \bS_{[-T,0]}(\mathsf q) + \bS_{[0,T]}(\mathsf q) + C_{\varrho} $$
uniformly in the path $\mathsf q$.
We use this in (\ref{gaussint3}), apply the Markov property of Brownian motion in the resulting term and reverse time in one of the factors to get
$$\norm[\cT]{P_T1}^2 \leq e^{C_\varrho} \left(\E_\W^0 (e^{\bS_{[0,T]}}) \right)^2 = e^{C_\varrho} \skalprod[\cT]{1}{P_T1}^2,$$
and thus 
\begin{equation} \label{not zero}
\skalprod[\cT]{1}{\frac{P_T1}{\norm[\cT]{P_T1}}} \geq e^{-C_{\varrho}/2}.
\end{equation}
The family $(P_T1/\norm[\cT]{P_T1})_{T>0}$ is bounded and thus relatively compact in the weak topology of $\cT$. Let $\Psi$ be the weak limit along a subsequence $(T_n)$ with $T_n \to \infty$ as $n \to \infty$. By (\ref{not zero}),
$\Psi \neq 0$. We want to establish that $\Psi$ is an eigenvalue of $P_{t}$. Let $\mu$ be the spectral measure of $H$ with respect to the vector $1 \in D(H)$ (domain in $\cT$). 
Then the infimum of the support of $\mu$ is  
$$E_{0} = \inf\mathrm{\supp}(\mu) = - \lim_{T \to \infty} \frac{1}{T} \ln \left(\int_{-\infty}^{\infty} e^{-Tx} \, d\mu(x) \right).$$
$E_0$ is also the infimum of the spectrum of $H$. This follows from the fact that $P_t$ preserves positivity: for $f \in L^{\infty} \cap D(H)$ we have
$$ \skalprod[\cT]{f}{P_{t}f} \leq \skalprod[\cT]{|f|}{P_t|f|} \leq \norm[L^{\infty}]{f}^2 \skalprod[\cT]{1}{P_t1}$$
for all $t > 0$, and consequently the infimum of the support of the spectral measure associated to $f$ is greater than $E_0$. Since $L^{\infty} \cap D(H)$ is dense, $E_0$ must be the infimum of the spectrum of $H$.
Again by using the spectral measure, we find
$$ \lim_{T \to \infty} \frac{\norm[\cT]{P_T1}}{\norm[\cT]{P_{T+s}1}} = e^{E_0s},$$
and thus $\skalprod[\cT]{\Psi}{P_t\Psi} = e^{-E_0t} = \norm{P_t}$. This implies that $\Psi$ is an eigenfunction of $P_t$, and thus of $H$. Since $P_t$ is positivity improving, $\Psi$ is unique and can be chosen strictly positive by the Perron-Frobenius theorem (cf. \cite{Reed/Simon}, section XIII.12, vol. 4).
\end{Proof}

Having existence and uniqueness of $\Psi$ under control, the spectral theorem yields
\begin{equation} \label{convergence in L^2}
\Psi_T = e^{TE_0}(P_T1) \to \skalprod[\cT]{1}{\Psi}\Psi \quad \mbox{as }T \to \infty
\end{equation}
in $\cT$ and thus in $L^\infty(\R^d,L^2(\g))$. 

In the following, $\Psi$ will be chosen strictly positive and normalized.
In the context of the Nelson model, $\Psi$ is the ground state of the dressed electron 
for total momentum zero. Its existence (also for small nonzero momentum) was first proven by Fr\"ohlich \cite{Froehlich} using a completely different method.

It is now easy to identify the infinite volume limit of the families $\P_{T}$ and $\N_{T}$. For an interval $I \subset \R$ let us denote by $\F_{I}$ the $\sigma$-field generated by the point evaluations 
$\{ (q_t, \phi_t): t \in I \}$.
Let $\P$ be the probability measure on paths $(q_t,\phi_t)_{t \in \R}$ determined by 
\begin{equation} \label{P}
\E_{\P}(f) = e^{2TE_0} \E_{\W \otimes \G} \left( \left. \Psi(q_{-T},\phi_{-T})
e^{- \int_{-T}^T \tau_{q_s}\phi_s( \varrho) \, ds} 
\Psi(q_{T},\phi_{T}) f \right| q_0 = 0 \right)
\end{equation}
for each bounded, $\F_{[-T,T]}$-measurable function $f$. From Theorems \ref{strongly continuous} and \ref{existence of ground state} we conclude that $\P$ is the measure of a Markov process with generator $L$ acting as
\begin{equation} \label{process generator}
Lf = - \frac{1}{\Psi}(H-E_{0})(\Psi f).
\end{equation}
This is in complete analogy with the ground state transform known from Schr\"odinger semigroups. 

\begin{Proposition} \label{P_T to P}
$\P_T \to \P$ as $T \to \infty$ in the topology of local convergence, i.e. $\E_{\P_T}(f) \to \E_\P(f)$ for each bounded, $\F_{[-t,t]}$-measurable function $f$ and each $t > 0$.
\end{Proposition}

\begin{Proof}
Defining $\Psi_T$ as in (\ref{convergence in L^2}), we have for $f \in \F_t$ and $t<T$
$$ \E_{\P_T}(f) = \frac{e^{2tE_0}}{\norm[\cT]{\Psi_T}^2} \E_{\W \otimes \G} \left( \left. \Psi_{(T-t)}(q_{-t},\phi_{-t})
e^{ - \int_{-t}^t \tau_{q_s} \phi_s(\varrho) \, ds} f \Psi_{(T-t)}(q_{t},\phi_{t})
\right| q_0=0 \right) $$
By (\ref{convergence in L^2}) and $\skalprod[\cT]{1}{\Psi} > 0$, $ \Psi_{(T-t)} / \norm[\cT]{\Psi_T} \to \Psi$ in $\cT$. Moreover, for bounded $f \in \F_{[0,T]}$, the map $Q$ on $L^{\infty}(\R^d,L^2(\g))$ with
$$ (Qg)(q, \phi) = \E_{\W \otimes \G}^{q,\phi} \left( \exp \left( - \int_{0}^t \tau_{q_s}\phi_s( \varrho) \, ds\right) f(q,\phi) 
g(q_t,\phi_t) \right) $$
is a bounded linear operator on $L^{\infty}(\R^d,L^2(\g))$; this follows from
$ |Qg| \leq \norm[\infty]{f} P_t|g| $ and the boundedness of $P_t$. Thus 
$$Q (\Psi_{(T-t)} / \norm[\cT]{\Psi_T}) \stackrel{T \to \infty}{\to} Q\Psi$$
in $L^{\infty}(\R^d,L^2(\g))$, and the claim follows.
\end{Proof}

\begin{Theorem} \label{integrability}
The family $(\N_T)_{T > 0}$ converges to a probability measure $\N$ in the topology of local convergence. Moreover, if $f \in \F_{[-t,t]}$ depends only on $\{q_s:-t\leq s \leq t\}$ and satisfies 
$\E_{\W}^0(|f|) < \infty$, then also 
$\E_\N(|f|) < \infty$, and $\E_{\N_T} (f) \to \E_\N(f)$ for such $f$.
\end{Theorem}

\begin{Proof}
The first statement follows from Proposition \ref{P_T to P} when considering functions of $q$ only. 
All the other statements will be proved once we show that there exists $C > 0$ such that 
\begin{equation} \label{uniform bound}
\sup_{T > 0} \E_{\N_T}(|f|) \leq C \E_\W^0(|f|)
\end{equation}
for  all $f \in \F_{[-t,t]}$. To see (\ref{uniform bound}), first note that
$$
 \E_{\N_T}(|f|) \leq  \frac{e^{2C_{\varrho}}}{Z_T} \E_\W^0  \left(e^{\bS_{[-T,-t]}}e^{\bS_{[-t,t]}}e^{\bS_{[t,T]}} |f| \right). 
 $$
Using the Markov property, stationarity of increments and time reversal invariance of two-sided Brownian motion, the latter term above is equal to 
$$e^{2C_{\varrho}} \E_\W^0 \left( R_T(q_{-t}) e^{\bS_{[-t,t]}(q)}f(q) R_T(q_t) \right),$$
where
 $$ R_T(q) = \frac{1}{\sqrt{Z_T}} \E_\W^q( e^{\bS_{[0,T-t]}}) = \frac
 {\skalprod[\cT]{1}{P_{T-t}1}}{\norm[\cT]{P_T1}}. $$
$R_T(q)$ is therefore independent of $q$ and convergent as $T \to \infty$. Since the pair potential $W$ is bounded, also $\bS_{[-t,t]}$ is uniformly bounded, and (\ref{uniform bound}) follows.
\end{Proof}

$\N$ describes the evolution of a particle driven by a stationary field. This field is given by the $\cK$-valued  process 
$$ \eta_t = \tau_{q_t} \phi_t.$$
$\eta_t$ is the configuration of the field $\phi_t$ as seen from the location $q_t$ of the particle. 

\begin{Proposition} \label{reversible}
The process $(\eta_t)_{t \in \R}$ is a reversible Markov process under $\P$. The reversible measure is given by $(U^{-1}\Psi)^2 \g$.
\end{Proposition}

\begin{Proof}
Let $f,g \in L^2(\g)$. Then (\ref{P}) implies
$$ \E_{\P}( f(\eta_{s}) g(\eta_{t})) = e^{|t-s|E_{0}}\skalprod[\cT]{\Psi Uf}{P_{|t-s|}(\Psi Ug)}.$$
Thus the generator of the $\eta_{t}$-process is unitarily equivalent to the operator $L$ (cf. (\ref{process generator})) on the Hilbert space $(\cT,\norm[\Psi]{.})$, where 
$\norm[\Psi]{f} =  \norm[\cT]{\Psi f}$. $L$ is self-adjoint on this Hilbert space, $L1 = 0$, and $\norm[\Psi]{1} = 1$. This proves reversibility.
\end{Proof}

The significance of the process $(\eta_t)_{t \in \R}$ is that it governs the infinitesimal increments of the process $(q_t)_{t \in \R}$ under $\N$. More explicitly, let  $\gamma \in \R^d$ be fixed, and 
$h_\gamma(q) = \gamma \cdot q$. From Proposition \ref{action of L} it will follow that $Lh_\gamma(q,\phi) = (\gamma \cdot \nabla_q \ln(\Psi))(q,\phi)$, and thus $Lh_\gamma \in \cT$. 
With 
\begin{equation}\label{Lf=j}
j = U^{-1}(\gamma \cdot \nabla_q \ln \Psi) \in L^2(\g),
\end{equation}
we have $L(\gamma \cdot q) = j(\eta)$. This fact is paraphrased by saying that $q_t$ is driven by $\eta_t$.

\section{The central limit theorem}
In this section we prove a functional central limit theorem for the process $(q_{t})_{t \geq 0}$ under the measure $\P$.  
The results of the previous section make it possible to apply the Kipnis-Varadhan method. The first ingredient is 

\begin{Theorem} \label{KVTheorem} \cite{KV} Let $(y_t)$ be a Markov 
process with respect to a filtration $\F_{t}$. Assume that $(y_{t})$ is reversible
with respect to a probability measure $\mu_0$, and that the
reversible stationary process $\mu$ with invariant measure $\mu_0$ is
ergodic. Let $V$ be a $\mu_0$ square integrable function on the state space
with $\int V \, d\mu_0 = 0$. Suppose in addition that $V$ is in the domain
of $L^{-1/2}$, where $L$ is the generator of the process $y_t$. Let 
$$ X_t = \int_0^t V(y_s) \, ds.$$ 
Then there exists a square integrable martingale $(N_t, \F_t)$ such that 
$$ \lim_{t \to \infty} \frac{1}{\sqrt t} \sup_{0 \leq s \leq t} |X_s - N_s| =
0$$ in probability with respect to $\mu$, where $X_0 = N_0 = 0$. Moreover, 
$$ \lim_{t \to \infty} \frac{1}{t} \E_{\mu} (|X_t - N_t|^2) = 0.$$
\end{Theorem}

In \cite{KV} the theorem is stated for the case that $\F_{t}$ is the filtration generated by the process $(y_{t})$. From the proof given there, it is obvious that the above modification of the theorem is also valid.

Let again $L$ be the generator of the $\P$-process, cf. (\ref{process generator}), and let $h_\gamma (q) = \gamma \cdot q$ with
$\gamma \in \R^d$ fixed. 
We write 
\begin{equation} \label{sum of two martingales}
 \gamma \cdot q_t = \left(\gamma \cdot q_{t} - \int_{0}^{t} Lh_\gamma(q_s, \phi_s) \, ds \right) +  \int_{0}^{t} Lh_\gamma(q_s, \phi_s)) \, ds.
 \end{equation}
The term in brackets will turn out to be an $\F_t$-martingale, and Theorem \ref{KVTheorem} will be applicable to the remaining term on the right hand side. This gives $\gamma \cdot q_t$ as the sum of two martingales, and the martingale central limit theorem may be applied. We now elaborate this program and start with a result that allows us to calculate $L(\gamma \cdot q)$.

\begin{Proposition} \label{action of L}
If $g \in C^2(\R^d)$, then $\Psi Lg \in C(\R^d, L^2(\g))$. Moreover,
$$
Lg(q,\phi) = \frac12 \Delta g(q) + \nabla_q g(q) \cdot \nabla_q \ln \Psi(q,\phi).
$$
\end{Proposition}

\begin{Proof}
For $\alpha \in \{1,\ldots,d\}$, $\partial_\alpha = \frac{\partial}{\partial q_\alpha}$ is the generator of the unitary group $f \mapsto \tau_{tq_\alpha}f$ on $\cT$, and is thus anti-selfadjoint on $\cT$. By Proposition \ref{kr}, $H_{\mr f} + V_{\varrho}$ is bounded below on $\T$ by $-a \in \R$, say. Thus 
\begin{eqnarray*}
E_0 & = & \skalprod[\cT]{\Psi}{H\Psi} = -\frac 12 \sum_{\alpha=1}^d \skalprod[\cT]{\Psi}{\partial_\alpha^2\Psi} + \skalprod[\cT]{\Psi}{(H_{\mr f} + V_{\varrho})\Psi} \geq \\
& \geq & \frac 12 \sum_{\alpha=1}^d \norm[\cT]{\partial_\alpha\Psi}^2  - a.
\end{eqnarray*}
This shows $\partial_\alpha \Psi \in \cT$. Now by $(H-E_0)\Psi = 0$ and $(H_{\mr f} + V_{\varrho})g\Psi = g (H_{\mr f} + V_{\varrho})\Psi$, we find 
$$ (H-E_0)g\Psi = - \frac 12  \Psi \Delta g - \nabla_q g \cdot \nabla_q \Psi,$$
and the proof is complete.
\end{Proof}

\begin{Proposition} \label{martingale1}
Let $j$ be defined as in (\ref{Lf=j}). Then 
\begin{equation} \label{martingale}
 M_t = \gamma \cdot q_t - \int_0^t j(\eta_s) \, ds 
\end{equation}
is an $\F_t$-martingale with stationary increments under $\P$. The quadratic variation of $M_t$ is $|\gamma|^2t$.
\end{Proposition}

\begin{Proof}
Let $f \in C(\R^d,\R)$ such that $E_\P(|f(q_t)|^2) < \infty$ for all $t \geq 0$, and suppose that $Lf$ and $L(f^2)$ are in all the $L^2$-spaces induced by the images of $\P$ under 
$(q, \phi) \mapsto (q_t,\phi_t), t \geq 0$. Then 
$$
\Pi_tf - f = \int_0^t \Pi_s Lf \, ds
$$
in all the $L^2$-spaces introduced above. In the terminology of \cite{Ethier/Kurtz}, $(f,Lf)$ is in the full generator of the transition semigroup $\Pi_t = e^{tE_0}\frac{1}{\Psi}P_t\Psi$ of the $\P$-process. It follows that $M_t = f(q_t) - \int_0^t Lf(q_s,\phi_s) \, ds$ is a square integrable martingale, and a direct calculation gives 
$$\E_{\P}((M_t)^2) = \E_\P(M_0^2) + \int_0^t \E_\P\left(L(f^2)(q_s,\phi_s) - 2 f(q_s) Lf(q_s,\phi_s)\right) \, ds.$$
Using this general theory and the fact that $L(\gamma \cdot q) = j(\eta)$, we just have to check that $f(q) = \gamma \cdot q$ fulfills the integrability conditions required at the beginning of the proof. Using Proposition \ref{action of L} and Theorem \ref{integrability}, this is immediate.
\end{Proof}

\begin{Proposition}
As $t \to \infty$,
\begin{equation} \label{diffusion constant}
 \lim_{t \to 0} \frac{1}{t} \E_\P((\gamma \cdot q_t)^2) = |\gamma|^2 - 2 \skalprod[\cT]{\gamma \cdot \nabla_q\Psi}{(H-E_0)^{-1}\gamma \cdot \nabla_q\Psi}.
 \end{equation}
\end{Proposition}
\begin{Proof}
By (\ref{martingale}), 
\begin{equation} \label{eqt}
\E_\P((\gamma \cdot q_t)^2) = \E_\P(M_t^2) - \E_\P\left(\left( \int_0^t j(\eta_s) \, ds \right)^2\right) + 2 \E_\P\left((\gamma \cdot q_t)
\int_0^t j(\eta_s) \, ds \right).
\end{equation}
The third term in (\ref{eqt}) is zero. This can be seen as follows: We have
$$ \E_\P\left((\gamma \cdot q_t) \int_0^t j(\eta_s) \, ds \right) = \E_\W^0(  (\gamma \cdot q_t) I(\mathsf q)),$$ 
where 
$$ I(\mathsf q) = \E_\G \left( \Psi(\eta_0)e^{-\int_0^t \eta_s(\varrho) \, ds} \left( \int_0^t j(\eta_s) \, ds \right) \Psi(\eta_t) \right)$$
and $\mathsf q$ denotes a path $(q_t)_{t \in \R}$.
Put $\tilde{q}_s = q_{t-s} - q_t$. Then by the reversibility of $\G$ and the fact that $\G$ is invariant under  the constant shift by $\tau_{q_t}$, we have $I(\tilde{\mathsf q}) = I(\mathsf q)$. Moreover $\tilde{q}_t = -q_t$, $\W^0$-almost surely, and $\W^0$ is invariant under the transformation $\mathsf q \mapsto \tilde{\mathsf q}$. Thus 
$$\E_\W^0(  (\gamma \cdot q_t)I(\mathsf q)) = - \E_\W^0( (\gamma \cdot q_t) I(\mathsf q) ) = 0.$$ 
From Proposition \ref{martingale1} we know that $\E_\P(M_t^2) = |\gamma|^2t$. 
Let $\Pi_t$ denote again the transition semigroup of $\P$. Then
\begin{eqnarray*}
 \frac{1}{t} \E_\P\left(\left( \int_0^t j(\eta_s) \, ds \right)^2\right) & = & \frac{1}{t} \int_0^t ds \int_0^t dr \, \skalprod[\Psi_0]{\gamma \cdot \nabla_q \ln \Psi}{\Pi_{|r-s|}\gamma \cdot \nabla_q \ln \Psi} \\
 & \stackrel{t \to \infty}{\to} & -2 \skalprod[\Psi_0]{\gamma \cdot \nabla_q \ln \Psi}{L^{-1}\gamma \cdot \nabla_q \ln \Psi} = \\
 & = &  2 \skalprod[\cT]{\gamma \cdot \nabla_q\Psi}{(H-E_0)^{-1}\gamma \cdot \nabla_q\Psi}.
 \end{eqnarray*}
Note that the last quantity is automatically finite; this follows from (\ref{eqt}) and the positivity of $P_t$. The proof is completed.
\end{Proof} 

\begin{Theorem} \label{clt}For the process $q_t$ under $\P$, a functional central limit theorem holds, i.e. as $\varepsilon \to 0$, the process $t \mapsto \sqrt{\varepsilon} q_{t/\varepsilon}, t \geq 0,$ converges in distribution to Brownian motion with diffusion matrix $D$ given by 
$$ D_{\alpha\beta} = \delta_{\alpha\beta} - 2 \skalprod[\cT]{\partial_\alpha\Psi}{(H-E_0)^{-1}\partial_\beta\Psi}, \quad (\alpha,\beta = 1,\ldots,d).$$
\end{Theorem}
\begin{Proof}
We have to check that the process $\eta_s$ and the functional $j$ fulfill the assumptions of Theorem \ref{KVTheorem}. If so, (\ref{sum of two martingales}) shows that $q_t$ is the sum of two martingales with stationary increments and one negligible process, and the martingale functional central limit theorem can be applied. The diffusion matrix can then be obtained from (\ref{diffusion constant}) by choosing the canonical basis vectors $e_\alpha$ in place of $\gamma$ and polarization. For checking conditions of Theorem \ref{KVTheorem}, note that $\E_\P(j(\eta_t)^2) < \infty$ was shown in the proof of Proposition \ref{action of L}, 
and (\ref{diffusion constant}) implies $j \in D(L^{-1/2})$. Moreover, $\E_\P(\eta_t) = \sum \gamma_\alpha \skalprod[\cT]{\Psi}{\partial_\alpha \Psi} = 0$ since $\partial_\alpha$ is anti-selfadjoint on $\cT$ and $\Psi$ is real-valued. Finally, $\eta_s$ is a reversible Markov process with respect to $\F_t$ by Proposition \ref{reversible}, and is ergodic, since $P_t$ is positivity improving.
\end{Proof}

\section{Lower bound for the diffusion matrix}

To complete the proof of Theorem \ref{a.the2} we still have to
show that $D>0$. For this purpose we follow ideas from Brascamp et
al. \cite{Brascamp}, who study fluctuations for anharmonic lattices, one
particular case of which is the path measure
$\mathcal{N}_T$ when discretized. 
Our main effort is to show that the desired lower
bound survives in the limits of zero discretization and of $T \to
\infty$.

Let us discretize $[-T,T]$ with step-size $\varepsilon$,
$N \varepsilon =T$, and approximate (\ref{a.j}) through
$$
\frac{1}{Z_N}
\exp[-H_\varepsilon (x)]  \prod^N_{{j=-N} \atop {j\not=0}}dx_j = \mu_N
$$
as a probability measure on $\mathbb{R}^{2Nd}$ with 
$$
H_\varepsilon(x) = \frac{1}{2\varepsilon} \sum^{N-1}_{j=-N}
(x_{j+1}-x_j)^2 + \frac{1}{2} \kappa\varepsilon\sum^N_{j=-N} x^2_j
+\frac{1}{2} \varepsilon^2 \sum^N_{i,j=-N}
W(x_i-x_j,\varepsilon(i-j))\,,
$$
$\kappa>0$. Here $x=(x_{-N},
\ldots, x_N)$, $x_0=0$. Clearly, $\mu_N \to \mathcal{N}_T$ weakly in the
limits $\kappa\to 0$, $\varepsilon \to 0$. Expectations with
respect to $\mu_N$ are denoted by $\E_N$.

We define the $2Nd \times 2Nd$ matrix $M^\kappa$ through
\begin{equation}\label{d.c}
M^\kappa_{i\alpha,j\beta} = \expval[N]{ \partial_{i\alpha} H
\partial_{j\beta} H} = \expval[N] {\partial_{i\alpha} \partial_{j\beta}
H } \,,
\end{equation}
$i,j = -N, \cdots, N$, $i,j \not=0$, $\alpha,\beta=1,\cdots, d$,
$\partial_{i\alpha} = \partial/\partial x_{i\alpha}$. Let
$\sum^\ast_{i,\alpha}$ denote the sum $i=-N, \cdots, N$, $i \not=
0$, $\alpha=1, \cdots, d$. Then for real coefficients
$f_{i\alpha}$, $g_{i\alpha}$ one has, by partial
integration and Schwarz inequality,
\begin{eqnarray*}
\Big( {\sum_{i,\alpha}}^\ast f_{i\alpha} g_{i\alpha}\Big)^2
= \E_N \Big( \Big( {\sum_{i,\alpha}}^{\ast} f_{i\alpha}
x_{i\alpha}\Big) \Big( {\sum_{j,\beta}}^\ast g_{j\beta}
\partial_{j\beta}H \Big) \Big)^2  \\
\leq \E_N \Big( \Big( {\sum_{i,\alpha}}^\ast f_{i\alpha}
x_{i\alpha}\Big)^2\Big)  \E_N \Big( \Big({\sum_{j,\beta}}^\ast g_{j\beta}
\partial_{j\beta}H\Big)^2\Big)\,.
\end{eqnarray*}
Since $M^{\kappa} \geq \kappa >0$, we can set
$g=(M^{\kappa})^{-1}f$ to obtain
\begin{equation}\label{d.e}
\E_N\Big({\Big(({\sum_{i,\alpha}}^\ast f_{i\alpha} x_{i\alpha}\Big)^2}\Big) 
\geq {\sum_{i\alpha,
j\beta}}^\ast \left((M^{\kappa})^{-1}\right)_{i\alpha,j\beta}
f_{i\alpha}f_{j\beta}\,.
\end{equation}

From (\ref{d.c}) one obtains $M^{\kappa}_{i\alpha,j\beta} =
M_{i\alpha,j\beta}+ \kappa \varepsilon \delta_{ij}
\delta_{\alpha\beta}$ with
$$
M_{i\alpha,j\beta}= -\varepsilon^{-1}\Delta^0_{ij}
\delta_{\alpha\beta} + \varepsilon^2 \Big( K^T_\varepsilon (i\alpha,
j\beta) - \delta_{ij} \sum^N_{n=-N\atop{n \neq 0}} K^T_\varepsilon
(i\alpha, n\beta) \Big)\,,
$$
where $\Delta^0$ is the lattice Laplacian with Dirichlet boundary
condition at $j=0$, Neumann boundary condition at $j=\pm N$, and
$$
K^T_\varepsilon(i\alpha,j\beta) = \expval[N]{\partial_\alpha
\partial_\beta W(x_i -x_j,\varepsilon(i-j)}\,,
$$
$\partial_\alpha = \partial/\partial x_\alpha$. In (\ref{d.e}) we
further decrease on the right by substituting $(M^\kappa +
\lambda)^{-1}$ for $(M^\kappa)^{-1}$, $\lambda>0$. We now take the
limit $\varepsilon \to 0$, $\kappa \to 0$. Then, for $f \in C
([-T,T], \mathbb{R}^d)$, we conclude that
\begin{equation}\label{d.h}
\mathbb{E}_{\N_T} \left(\left(\int^T_{-T} f(t) \cdot q_t \, dt\right)^2\right) \geq \left \langle
f, (A^0_T + \lambda)^{-1} f \right \rangle
\end{equation}
for all $\lambda>0$. Here $\langle\cdot,\cdot\rangle$ denotes the
inner product for $L^2([-T,T],dt)\otimes \mathbb{C}^d$, $A^0_T$
is the linear operator
$$
A^0_T = -\Delta^{0,T}\otimes 1 + B^T\,,
$$
$-\Delta^{0,T} = -d^2/dt^2$ with Dirichlet boundary condition at
0 and Neumann boundary condition at $\pm T$, and $B^T$ is the
integral operator
$$
B^T f_\alpha (t) = \sum^d_{\beta=1} \int^T_{-T}
K^T_{\alpha\beta}(t,s)(f_\alpha(t) - f_\beta(s))ds
$$
with kernel
$$
K^T_{\alpha\beta}(t,s) = \mathbb{E}_{\N_T} (\partial_\alpha
\partial_\beta W (q_t -q_s, t-s))\,.
$$

We take the limit $T\to\infty$ in (\ref{d.h}). For compactly supported $f$, 
the left hand
side converges to $\mathbb{E}_\N ((\int dtf(t) \cdot q_t)^2)$  by Theorem \ref{integrability}. For
the right hand side we use a theorem by Kurtz (\cite{Ethier/Kurtz}, Thm. 1.6.1). 
We consider $L^2 ([-T,T]) \otimes \mathbb{C}^d$ as a subspace on $L^2
(\mathbb{R}) \otimes \mathbb{C}^d$. 
Let $A^0 = -\Delta^0 \otimes 1 + B$
where $-\Delta^0$ is the Laplacian on $L^2(\mathbb{R})$ with
Dirichlet boundary conditions at $t=0$ and where
$$
B f_\alpha (t) = \sum^d_{\beta=1} \int
K_{\alpha\beta}(t,s)(f_\alpha(t) - f_\beta(s))ds
$$
with kernel
$$
K_{\alpha\beta}(t,s) = \mathbb{E}_\N (\partial_\alpha
\partial_\beta W (q_t -q_s, t-s))\,.
$$
Putting $K_{\alpha\beta} = K^\infty_{\alpha\beta}$, we find for all $T \leq \infty$ 
$$
 \left| \int K_{\alpha \beta}^T(t,s) \, ds \right| \leq \int \frac{\widehat{\varrho}(k)^2}{2\om(k)^2}|k|^2 \, dk
$$
and 
$$
 \left| \int K_{\alpha \beta}^T(t,s) f_\beta(s) \, ds \right| \leq \int  \int \frac{\widehat{\varrho}(k)^2}{2\om(k)}|k|^2 e^{-\om(k)|t-s|} |f_{\beta}(s)| \, dk \, ds.
$$
 By the Schwarz inequality, the $L^2(\R)$-norm of the right hand side above is bounded by $\norm[L^2(\R)]{f} \int \frac{\widehat{\varrho}(k)^2}{2\om(k)^2}|k|^2 \, dk$, and thus $B_T$ and $B$ are bounded operators on $L^2(\R) \otimes \C^d$. Since clearly $\lim_{T\to\infty}|Bf_\alpha(t) - B^Tf_\alpha(t)| = 0$ pointwise, dominated convergence yields $\lim_{T\to\infty}B^Tf = Bf$ in $L^2(\R^d) \otimes \C$ for all $f \in L^2(\R^d) \otimes \C$. Now let
$$D=\{g
\in L^2(\mathbb{R})\otimes \mathbb{C}^3|
g_\alpha(0)=0,g''_\alpha \in L^2; g_\alpha \mbox{ has compact
support, } \alpha=1, \cdots,d\}.$$ 
$D$ is a core for $A^0$, and $\lim_{T\to \infty}
\norm{((-\Delta^{0,T} + \Delta^0)\otimes 1) g}=0$. Thus also $\lim_{T \to \infty} A^0_Tg = A^0g$, and from (\ref{d.h}) and  \cite{Ethier/Kurtz}, Thm. 1.6.1 we conclude
\begin{equation}\label{d.p}
\mathbb{E}_\N \Big( \Big( \int f(t) \cdot q_t \, dt \Big)^2 \Big) \geq \langle f, (A^0
+ \lambda)^{-1}f \rangle
\end{equation}
for $\lambda>0$ and $f \in L^2(\mathbb{R}) \otimes \mathbb{C}^d$.

It remains to estimate the right hand side of (\ref{d.p}).
Let $\widehat K_{\alpha\beta}$ be the Fourier transform of $K$.
$K$ only depends on $t-s$, since $\N$ has stationary
increments, thus in Fourier space $B$ is multiplication by $\widehat{K}(0) -
\widehat{K}(\omega)$. Since by assumption (\ref{cond3}) $\int
t^2|K_{\alpha\beta} (t)|dt < \infty$, there exists a constant
$D_0$ such that
\begin{equation}\label{d.q}
\widehat K(0)-\widehat K(\omega) \leq \frac {1}{1+\omega^2} D_0
\omega^2
\end{equation}
as $d \times d$ matrices. Let $\bar{B}$ be the operator corresponding to the right hand side of (\ref{d.q})
 and $-\Delta^0 + \bar{B}=\bar{A}^0$. By operator
monotonicity and taking $f_\alpha \to \gamma \delta_t$ one achieves
\begin{equation}\label{d.r}
\mathbb{E}_\N \Big( (\gamma \cdot q_t)^2 \Big) \geq  |\gamma|^2 (\bar{A}^0 + \lambda)^{-1}
(t,t)\,.
\end{equation}

Let $A= -\Delta + \bar{B}$ on $L^2(\mathbb{R})$. Then, as can be
checked directly,
$$
(\bar{A}^0 + \lambda)^{-1} (t,s) = (A + \lambda)^{-1} (t,s)-\left[ (A +
\lambda)^{-1}) (0,0)\right]^{-1} (A + \lambda)^{-1} (0,s) (A +
\lambda)^{-1} (t,0)\,.
$$
$A$ is multiplication by $\omega^2+(1+\omega^2)^{-1} D_0 \omega^2$
in Fourier space. Thus by explicit computation
\begin{equation}\label{d.t}
\lim_{\lambda \to 0}(\bar{A}^0 + \lambda)(t,t) =
(1+D_0)^{-1}|t|\,.
\end{equation}

Combining (\ref{d.r}) and (\ref{d.t}) means that there is a
constant $c_0$ such that $c_0>0$ and
$$
\mathbb{E}_\N\Big( (\gamma \cdot q_t)^2 \Big) \geq c_0 |\gamma|^2 |t|\,
$$
for all $\gamma \in \R^d$. Hence $D\geq c_0 >0$.

\section{Appendix: Proof of Theorem \ref{strongly continuous}}

We will prove Theorem \ref{strongly continuous} and start with collecting some facts that will be useful later on. 
\begin{Proposition} \label{kr}
Let $H_{\mr f}$ and $V_{\varrho}$ be as in (\ref{generator}). For each $\varepsilon > 0$ there exists a constant $C$ independent of $q$ such that 
\begin{equation} \label{Kato-Rellich estimate}
\norm[L^2(\g)]{V_{\varrho} f}^2 \leq \varepsilon \norm[L^2(\g)]{H_{\mr f}f}^2 + C\norm[L^2(\g)]{f}^2
\end{equation}
for all $f \in D(H_{\mr f})$. $H_{\mr f} + V_{\varrho}(q)$ is self-adjoint on $D(H_{\mr f})$ and bounded below uniformly in $q$.
\end{Proposition}
The proof of (\ref{Kato-Rellich estimate}) is standard and is most conveniently done in the Fock space representation of $L^2(\g)$. We refer to \cite{thesis}. From (\ref{Kato-Rellich estimate}), the other assertions follow by the Kato-Rellich theorem \cite{Reed/Simon}.

Apart from the spaces introduced in the paragraph above Theorem \ref{strongly continuous}, we will need
\begin{eqnarray*}
L^{\infty}(\R^d,L^2(\g)) & = & \{ f:\R^d \times \cK \to \C: \mathrm{ess}\sup_{q \in \R^d} \norm[L^2(\g)]{f(q,.)} < \infty\},  \\
L^{\infty}(\R^d \times \cK) & = & \{ f:\R^d \times \cK \to \C: \mathrm{ess}\sup_{q \in \R^d, \phi \in \cK} f(q,\phi) < \infty\}, \quad \mbox{and} \\
L^2(\R^d,L^2(\g)) & = & \{ f:\R^d \times \cK \to \C: \int |f(q,\phi)|^2 \, dq \, d\g(\phi) < \infty \}.
\end{eqnarray*}
Obviously, $L^{\infty}(\R^d \times \cK) \subset L^\infty(\R^d,L^2(\g))$.

We now give two identities connected with the semigroup $P_t$ given in (\ref{semigroup equation}). By the Markov property and time reversibility of $\W$ and $\G$ we find 
\begin{equation} \label{e1}
 \norm[L^2(\g)]{P_tf}^2(q) = \E_{\W \otimes \G}\left( \left.  f(q_{-t}, \phi_{-t}) e^{\int_{-t}^t \tau_{q_s}\phi_s(\varrho) \, ds}  f(q_t,\phi_t) \right| q_0 = q \right)
 \end{equation}
for $f \in L^{\infty}(\R^d \times \cK)$. 
As a semigroup of operators on $L^2(\R^d,L^2(\g))$,  $P_t$ is self-adjoint and strongly continuous, and
\begin{equation} \label{fkn}
 P_tf(q,\phi) = e^{-tH}f(q,\phi) \quad \mbox{in } L^2(\R^d,L^2(\g)),
 \end{equation}
where $H$ is given by (\ref{generator}). In other words, $-H$ is the generator of $P_t$ in $L^2(\R^d,L^2(\g))$. (\ref{fkn}) is the Feynman-Kac-Nelson formula, and a proof can be given e.g. via the Trotter product formula, cf. \cite{thesis}. Note that the content of Theorem \ref{strongly continuous} is to extend (\ref{fkn}) to the spaces $C_0(\R^d,L^2(\g))$ and $\cT$ carrying the $\sup$-norm given in (\ref{sup norm}). 

\begin{Proposition} \label{semigroup}
\begin{itemize}
\item[a) ] $P_t$ is a semigroup of bounded operators on $L^{\infty}(\R^d, L^2(\g))$.
\item[b) ] If $f \in L^{\infty}(\R^d \times \cK) \cap C(\R^d,L^2(\g))$, then  $P_tf \in  C_b(\R^d,L^2(\g))$.
\item[c) ] If $f \in L^{\infty}(\R^d \times \cK) \cap C(\R^d,L^2(\g))$,then
$$  \norm[L^2(\g)]{P_tf - f} \stackrel{t \to 0}{\to} 0$$
uniformly on compact subsets of $\R^d$.
\end{itemize}
\end{Proposition}
\begin{Proof} a)
We only have to show boundedness of $P_t$, the semigroup property then follows from the Markov property of $\G$ and $\W$. Since $|P_tf| \leq P_t|f|$ pointwise, it is sufficient to consider positive functions. 
At first let $f \in L^{\infty}(\R^d \times \cK)$. 

Let us fix a path $\mathsf q: [-t,t] \to \R^d$, not necessarily continuous,
and define the bilinear form 
$$
\skalprod[\mathsf q]{f}{g} = \int f(q_{-t}, \phi_{-t}) e^{- \int_{-t}^{t} \tau_{q_s}\phi_s(\varrho) \, ds} g(q_t,\phi_t) \, d\G(\phi).
$$
For two paths $\mathsf q, \tilde{\mathsf q}$, 
\begin{equation} \label{quad form continuous}
 |\skalprod[\mathsf q]{f}{g} - \skalprod[\tilde {\mathsf q}]{f}{g}| \leq \norm[\infty]{f} \norm[\infty]{g} 
\E_\G \left( \left( e^{- \int_{-t}^{t} \tau_{q_s}\phi_s(\varrho) \, ds} - e^{- \int_{-t}^{t} \tau_{\tilde{q}_s}\phi_s(\varrho) \, ds}\right)^2 \right),
\end{equation}
and an explicit Gaussian integration shows that $\mathsf q \mapsto \skalprod[\mathsf q]{f}{g}$ is continuous from $L^\infty([-t,t],\R)$ to $\R$. Now let $\mathsf q$ be continuous, and define  
$$ q^{(n)}_s = \sum_{j=-n}^{n - 1} q_{tj/n}1_{[tj/n,t(j+1)/n[}(s).$$
Then 
\begin{equation} \label{product}
\skalprod[\mathsf q^{(n)}]{f}{g} = \skalprod[L^2(\g)]{f}{\prod_{j=-n}^{n-1} S_{t/n}^{q_{tj/n}}g}
\end{equation}
 with
$$
 S_t^qg (\phi) =  \E_\G^\phi \left( e^{- \int_{0}^{t} \tau_{q}\phi_s(\varrho) \, ds} g(\phi_t) \right)
 $$
By Proposition \ref{kr},  $e^{-t(H_{\mr f} + V_{\varrho}(q))}$ is a bounded operator on $L^2(\g)$ with $\norm{e^{-t(H_{\mr f} + V_{\varrho}(q))}} \leq e^{t\alpha}$ for some $\alpha \in \R$.
On the other hand,
$$ S_t^qg (\phi) =  e^{-t(H_{\mr f} + V_{\varrho}(q))}g(\phi)$$
in $L^2(\g)$. This is a variant of (\ref{fkn}), and can be proven in the same way. 
 Now by the Schwarz inequality and repeated use of the operator norm inequality, (\ref{product}) yields
$$
|\skalprod[\mathsf q^{(n)}]{f}{g}| \leq e^{2t\alpha} \norm[L^2(\g)]{f}(q_{-t})\norm[L^2(\g)]{g}(q_t)
$$
for bounded $f$ and $g$. Since $\mathsf q^{(n)} \to \mathsf q$ uniformly on $[-t,t]$ as $n \to \infty$, this remains valid when we replace $\mathsf q^{(n)}$ by $\mathsf q$. An application of monotone convergence now gives 
\begin{equation} \label{quad form bounded}
|\skalprod[\mathsf q]{f}{g}| \leq e^{2t\alpha} \norm[L^2(\g)]{f}(q_{-t})\norm[L^2(\g)]{g}(q_t)
\end{equation}
for all $q \in C(\R, \R^d), f,g \in L^2(\g)$. 
We use this in (\ref{e1}) and have established our first claim.\\
b) Let $q,r \in \R^d$ and $f,g \in L^{\infty}(\R^d \times \cK)$, and write $f_r(q,\phi) = f(q+r,\phi)$. Then
\begin{eqnarray*}
\lefteqn{\norm[L^2(\g)]{P_tf(q,.) - P_tf(q+r,.)}^2 = } \\
 &=& \E_\g \left( \left( \E_{\W \otimes \G}^{q,\phi} \left( e^{- \int_0^t \tau_{q_s}\phi_s(\varrho) \, ds}f(q_t,\phi_t) - e^{- \int_0^t \tau_{(q_s+r)}\phi_s(\varrho) \, ds}f_r(q_t,\phi_t) \right) \right)^2 \right) = \\
 & = & \E_\W^q \left( \skalprod[\mathsf q]{f}{f} 
  +  \skalprod[\mathsf q + r]{f_r}{f_r} - 2 \skalprod[\mathsf q + 1_{\{t\geq0\}}r]{f}{f_r} \right).
  \end{eqnarray*}
  For each path $\mathsf q$, each of the terms in the last line above converges to $\skalprod[\mathsf q]{f}{f}$ as $r \to 0$ by (\ref{quad form continuous}),  (\ref{quad form bounded}) and the continuity of $q \mapsto f(q,.)$. Thus the integrand converges to zero pathwise, and the second claim follows by dominated convergence. \\
c) Let
$$ Q_tf(q,\phi) = \E_{\W \otimes \G}^{q,\phi}( f(q_t,\phi_t)).$$
Then for bounded $f$, 
$$ \norm[L^{\infty}(\R^d,L^2(\g))]{Q_tf -P_tf} \stackrel{t \to 0}{\to} 0,$$
as a calculation similar to the one in (\ref{quad form continuous}) shows. Therefore we must only show that  $\norm[L^2(\g)]{Q_tf - f}$ vanishes uniformly on compact sets as $t \to 0$. Let us write $f_t(q,\phi) = f(q,\phi_t)$. 
By reversibility of $\W$ and $\G$ and the Cauchy-Schwarz inequality, 
$$ \norm[L^2(\g)]{Q_tf(q,.) - f(q,.)}^2 \leq \E_\W^q \left( \norm[L^2(\G)]{f_t(q_t,.) - f_0(q_0,.)}^2 \right).$$
Moreover, 
\begin{eqnarray*}
\norm[L^2(\G)]{f_t(q_t,.) - f_0(q_0,.)} & \leq & \norm[L^2(\G)]{f_t(q_t,.) - f_t(q_0,.)} + \norm[L^2(\G)]{f_t(q_0,.) - f_0(q_0,.)} = \\
& = & \norm[L^2(\g)]{f(q_t,.) - f(q_0,.)} + \norm[L^2(\G)]{f_t(q_0,.) - f_0(q_0,.)}.
\end{eqnarray*}
Thus it is enough to show that
\begin{equation} \label{strong continuity 1}
\norm[L^2(\G)]{f_t(q,.) - f_0(q,.)} \stackrel{t \to 0}{\to} 0 \quad \mbox{uniformly on compacts},
\end{equation}
and
\begin{equation} \label{strong continuity 2}
\E_\W^q \left(  \norm[L^2(\g)]{f(q_t,.)-f(q_0,.)}^i \right) \stackrel{t \to 0}{\to} 0, \quad \mbox{uniformly on compacts}, i=1,2.
\end{equation} 
In order to prove (\ref{strong continuity 1}), suppose there exist bounded sequences $(q_n) \subset \R^d$, $(t_n) \subset \R^+$ with $t_n \to 0$, and $\norm[L^2(\G)]{f_{t_n}(q_n,.) - f_0(q_n,.)} > \delta $ for all $n$. We may assume that $q_n$ converges to $q \in \R^d$. Then, with the notation introduced above, 
\begin{eqnarray*}
\delta &<& \norm[L^2(\G)]{f_{t_n}(q_n,.) - f_0(q_n,.)} \leq  \norm[L^2(\G)]{f_{t_n}(q_n,.) - f_{t_n}(q,.)}  + \\ 
&& + \norm[L^2(\G)]{f_{t_n}(q,.) - f_0(q,.)}  +  \norm[L^2(\G)]{f_{0}(q,.) - f_0(q_n,.)} = \\
& = &   \norm[L^2(\G)]{f_{t_n}(q,.) - f_0(q,.)} + 2  \norm[L^2(\g)]{f(q_n,.) - f(q,.)}.
\end{eqnarray*}
By choosing $n_0$ so large that $\norm[L^2(\g)]{f(q_n,.) - f(q,.)} < \delta/3$ for all $n > n_0$, we find that $\norm[L^2(\G)]{f_{t_n}(q,.) - f_0(q,.)} > \delta/3$ for all these $n$. However, $\norm[L^2(\G)]{f_{t_n}(q,.) - f_0(q,.)}$ must converge to zero by the strong continuity of the semigroup corresponding to $\G$, thus we have a contradiction. To prove (\ref{strong continuity 2}), note that
 $q \mapsto f(q,.)$ is uniformly continuous on compact sets. Thus for $\varepsilon > 0$ we may choose $\delta > 0$ such that $\norm[L^2(\g)]{f(q,.) - f(\tilde q,.)} < \varepsilon/2$ for $|q - \tilde q| < \delta$. Now by the properties of Brownian motion, there exists $t_0 > 0$ such that for $0 \leq t < t_0$, $\W^q(|q_t-q| > \delta) < \varepsilon / (2^{i+1} \norm[L^\infty]{f}^i)$. Then 
$\E_\W^q (\norm[L^2(\g)]{f(q_t,.)-f(q_0,.)}^i) < \varepsilon$ for such $t$ and uniformly in $q$, proving (\ref{strong continuity 2}).
\end{Proof}

\begin{Proposition}
$P_t$ is a strongly continuous semigroup of bounded operators on $C_0 = C_0(\R,L^2(\g))$.
\end{Proposition}
\begin{Proof}
We want to use Proposition \ref{semigroup} b) and c) for approximation and therefore must first show that $L^{\infty}(\R^d \times \cK)$ is dense in $C_0$. Let $f \in C_0$ and take $f^R = (f \wedge -R) \vee R$ for $R \geq 0$. 
Since $L^{\infty}(\cK)$ is dense in $L^2(\g)$, for each $q \in \R^d$ and each $\varepsilon > 0$,
\begin{equation} \label{pointwise finite}
 R_q(\varepsilon) = \inf \{ R \geq 0: \norm[L^2(\g)]{f(q,.) - f^R(q,.)} \leq \varepsilon \} < \infty.
 \end{equation} 
$R_q(\varepsilon)$ is bounded on compact subsets of $\R^d$, for otherwise we would find $(q_n) \subset \R^d$ with $q_n \to q$ and $R_n = R_{q_n}(\varepsilon) > n$ for all $n$. Choosing $n_0$ so large that $\norm[L^2(\g)]{f(q_n,.) - f(q,.)} < \varepsilon/3$ for all $n > n_0$, 
we would have
\begin{eqnarray*}
\varepsilon & = & \norm[L^2(\g)]{f(q_n,.) - f^{R_n}(q_n,.)} \leq \\
 & \leq & \norm[L^2(\g)]{f(q_n,.) - f(q,.)} + \norm[L^2(\g)]{f(q,.) - f^{R_n}(q,.)} + \norm[L^2(\g)]{f^{R_n}(q,.) - f^{R_n}(q_n.,)} \leq \\
 & \leq & \norm[L^2(\g)]{f(q,.) - f^{R_n}(q,.)} + 2\varepsilon/3 \leq \norm[L^2(\g)]{f(q,.) - f^{n}(q,.)} + 2\varepsilon/3
 \end{eqnarray*}
 for each $n > n_0$. This implies $R_q(\varepsilon/3) = \infty$, in contradiction to (\ref{pointwise finite}). Thus $R_q(\varepsilon)$ is bounded on compacts. However, since $f \in C_0$, $R_q(\varepsilon) = 0$ for $|q|$ large enough, and thus $R_q(\varepsilon)$ is bounded on all $\R^d$. Thus bounded functions are dense in $C_0$. 
 
Let us now show that $P_t$ leaves $C_0$ invariant. From (\ref{e1}) and (\ref{quad form bounded}) we see that
$$\norm[L^2(\g)]{P_tf(q,.)} \leq e^{Ct} \E_\W^q \left( \norm[L^2(\g)]{f(q_{-t},.)} \norm[L^2(\g)]{f(q_{t},.)} \right).$$
$\W^0(|q_t| \geq R)$ decays exponentially in $R$ for all $t$, and so does $\norm[L^2(\g)]{P_tf(q,.)}$. From Proposition \ref{semigroup} b) it follows that $P_tf \in C_b(\R^d,L^2(\g))$, and thus   $P_tf \in C_0$ for $f \in C_0 \cap L^\infty(\R^d \times \cK)$. By Proposition \ref{semigroup} a) and approximation, we obtain $P_t C_0 \subset C_0$. In a similar way, we obtain strong continuity from Proposition \ref{semigroup} c).
\end{Proof}

\begin{Proposition}
$P_t$ is a strongly continuous semigroup of bounded operators on $\cT$.
\end{Proposition}
\begin{Proof}
Let $f \in \cT$. By the invariance of $\G$ under the map $\phi \mapsto \tau_q \phi$ for each $q \in \R^d$ and the translation invariance of Wiener measure,
\begin{eqnarray*}
 (P_t f)(q,\phi) & = & \E_{\W \otimes \G}^{0, \tau_q \phi} \left( e^{-\int_0^t \tau_{q_s}\phi_s(\varrho) \, ds} \tau_{q_t} U^{-1} f(\phi_t) \right) = \\
 & = & U \left( \E_{\W \otimes \G}^{0, \phi} \left( e^{-\int_0^t \tau_{q_s}\phi_s(\varrho) \, ds} \tau_{q_t} U^{-1} f(\phi_t) \right) \right)
 \end{eqnarray*}
and thus $P_tf \in \cT$. It is easy to see that $L^{\infty} \cap \cT$ is dense in $\cT$. Strong continuity then follows from Proposition \ref{semigroup} c) and (\ref{scalar product in T}). 
\end{Proof}

It remains to show that the generator of $P_t$ is given by $-H$, with $H$ from (\ref{generator}). Write $- H_{C_0}$ for the generator in $C_0$ and $- H_\T$ for the generator in $\T$. If $f \in D(H_\T)$ or $f \in D(H_{C_0})$, then $fg \in D(H_{C_0})$ for each $g \in C^2(\R^d,\C)$ with compact support, and $fg \in D(H_{L^2})$, where $H_{L^2}$ is the generator of $P_t$ as a semigroup on $L^2(\R^d,L^2(\g))$. By the Feynman-Kac-Nelson formula (\ref{fkn}), we have $H_{L^2}fg = Hfg$ almost everywhere, and thus also $H_{C_0}fg = -Hfg$ in $L^{\infty}(\R^d,L^2(\g))$. Since $H_\T$, $H_{C_0}$ and $H$ are local in $q$, we may now use a smooth partition of unity to conclude $H_{\T}f = H_{C_0}f = Hf$. $H$ is obviously symmetric in $\T$, and is thus self-adjoint as the generator of a strongly continuous semigroup.
The proof of Theorem \ref{strongly continuous} is completed.
\hspace*{\fill} \ensuremath{\square}

\end{document}